\documentclass[review]{elsarticle}    % Enable this line and disable the
\usepackage{graphicx} % for pdf, bitmapped graphics files
\usepackage{amsmath} % assumes amsmath package installed
\usepackage{amssymb}  % assumes amsmath package installed
\usepackage{amsfonts}
\usepackage{amsthm}
\usepackage{color}
\usepackage{psfrag}
\usepackage{eurosym}
\usepackage{booktabs}
\usepackage{multirow}
\usepackage{pdflscape}
\usepackage{enumitem}
\biboptions{numbers,sort&compress} % per raggruppare le citazioni
\usepackage{longtable}
\usepackage{algorithmic}
\usepackage{textcomp}
\usepackage{xcolor}
\newtheorem{problem}{Problem}
\def\cU{\mathcal{U}}
\def\cR{\mathcal{R}}
\newcommand{\ds}{\displaystyle}
\usepackage{subcaption} % For subfigures
\usepackage[utf8x]{inputenc}
\usepackage{multicol}
\usepackage{mathtools}
\usepackage{multirow}
\usepackage{a4wide}  % ******* aggiunta per diminuire i margini *******

\newtheorem{remark}{Remark}
\newtheorem{example}{Example}

\newcommand{\ba}{\begin{array}}
\newcommand{\ea}{\end{array}}

\usepackage{color}  % per poter vedere le parti tagliate in blu
\newcommand{\T}{\mathcal{T}}

\newcommand{\I}{\mathcal{I}}

\newcommand{\U}{\mathcal{U}}

\graphicspath{{./Figs/}}
\journal{Energy}

\begin{document}
\begin{frontmatter}
\title{Optimal Operation of Renewable Energy Communities \\ under Demand Response Programs}
\author{Gianni Bianchini}\ead{giannibi@diism.unisi.it}    % Add the
\author{Marco Casini}\ead{casini@diism.unisi.it}              % e-mail address
\author{Milad Gholami}\ead{gholami@diism.unisi.it}
\address{Dipartimento di Ingegneria dell'Informazione e Scienze Matematiche\\Universit\`a di Siena, Via Roma 56, 53100 Siena, Italy}  % Please supply
\begin{keyword}                           % Five to ten keywords,
Renewable Energy Communities;
Demand Response;
Electrical Energy Storage;
Optimization 
\end{keyword}
\begin{abstract}
Within  the context of renewable energy communities, this paper focuses on optimal operation of producers equipped with energy storage systems in the presence of demand response. A novel strategy for optimal scheduling of the storage systems of the community members under price-volume demand response programs, is devised. The underlying optimization problem is designed as a low-complexity mixed-integer linear program that scales well with the community size. Once the optimal solution is found, an algorithm for distributing the demand response rewards is introduced in order to guarantee fairness among participants. The proposed approach ensures increased benefits for producers joining a community compared to standalone operation.
\end{abstract}
\end{frontmatter}

\section{Introduction} \label{sec:introduction}

To mitigate the environmental impact of electric energy systems and facilitate the transition toward net-zero CO$_2$ emissions \cite{greenDeal}, several paradigms and strategies have been proposed in the literature. Among these, renewable energy communities (RECs) and demand response (DR) programs stand out as particularly promising tools.

A REC consists of a collective of entities that engage in energy exchange through the power grid \cite{caramizaru2020energy,gjorgievski2021social}. The primary objective of a REC is to deliver social, economic, and environmental benefits to its members by effectively coordinating their load, generation, and storage assets \cite{EUDir944}. While national regulations often govern REC organization \cite{WHK20}, individual communities retain some autonomy as far as the operation strategy for their activities is concerned. One of the main objectives of such strategy is to optimize the overall welfare coming from active participation of the entities in the REC. In this respect, suitable rules must be implemented for the redistribution of benefits
among participants, ensuring that no member may find it advantageous to exit the community \cite{rogerio23}. Moreover, the redistribution policy should ensure fairness among participants \cite{fioriti21}.

Demand response programs represent another avenue for reducing environmental impact \cite{albadi2008summary,siano2014demand,Losi2015,hussain2015review,Mohammad2019}. In this framework, participants voluntarily adjust their load profiles in response to specific requests by an aggregator (e.g., the Distribution System Operator (DSO)), thereby providing ancillary services to the grid, such as peak power reduction and enhanced network stability \cite{ma2013demand}. In exchange, a monetary reward is granted to participants that fulfill a given DR request. In addition to the flexibility provided by load profile shaping, the optimal operation of electrical energy storage (EES) systems can contribute to achieving DR goals \cite{ma2016demand,zhang2017demand}. In this respect, RECs have been found to represent an important source of DR flexibility, especially when equipped with EES facilities \cite{mancarella19}. Optimized scheduling of a smart community with shared storage in the presence of DR is addressed in \cite{hou24}. A comprehensive review of DR programs and their potential application to RECs is available in the literature \cite{honarmand2021overview}.

This paper deals with the integration of incentive-based DR \cite{khajavi2011role,mohandes2020incentive} in the design of REC operation policies. In particular, we focus on the so-called price-volume model. In this paradigm, customers receive monetary incentives for maintaining consumption below a specified threshold during fixed time intervals. Price-volume DR has been successfully applied in various contexts, including load forecasting \cite{ruiz2015residential,gorria2013forecasting}, smart buildings \cite{bianchini2016demand,bianchini2019integrated}, and electric vehicle charging station management \cite{zanvettor2024AE}. Our study diverges from existing research which mainly focuses on load profile shaping, both in the general smart grid \cite{shao2011demand,jiang2014load,law2012demand} and in the specific REC context \cite{sangare23,pelekis23,hussain24}, by investigating how DR requests can be effectively met through optimized operation of individual EES systems associated with renewable generation sources (e.g., photovoltaic (PV) panels or wind turbines) in a REC. Specifically, we assume that the DSO communicates a price-volume DR program to the REC at the beginning of each day and grants an incentive if the overall community load/generation lies within prescribed thresholds in given time intervals. These thresholds align with periods where net load reduction supports grid stability. To the best of our knowledge, this problem remains unexplored in the literature.

\subsection{Paper contribution}
This paper addresses the design of an optimal scheduling policy enabling renewable energy producers equipped with battery energy storage systems (BESS) to participate in a price-volume DR program within a REC.
The main contributions of the work can be summarized in the following points.
\begin{itemize}
\item A novel formulation is introduced for the problem of optimal management of storage systems connected to renewable generation plants within a REC in the presence of price-volume DR requests. A specific feature of the proposed approach is that it always guarantees an additional economic benefit to producers participating in the REC with respect to autonomous (standalone) operation. The optimization problem stemming from such formulation turns out to be a low-complexity mixed-integer linear program (MILP) involving few binary variables, thus making the approach computationally feasible even for large REC memberships.

\item For the problem introduced above, two different objective functions are proposed and compared. One is oriented toward the maximization of the total community revenue from DR, while the other aims to maximize the total profit of REC producers.

\item An algorithm to partition the overall community DR reward among REC members is designed. Such an assignment guarantees distribution of the reward among participant entities according to a fairness principle.
\end{itemize}
Performance evaluation is carried out via extensive numerical simulations, showing pros and cons of the proposed approach and of the considered objective functions.

\subsection{Paper structure}
The paper is organized as follows: in Section~\ref{sec:modeling}, the considered problem is formulated and models of REC entities and DR requests are reported. In Section~\ref{sec:optimization}, the optimization problem that yields the storage system operation policy is designed, and an algorithm for providing a fair redistribution of the community DR reward is proposed. Numerical examples are presented in Section~\ref{sec:test_cases} and discussed in Section~\ref{sec:discussion}. Finally, conclusions are drawn in Section~\ref{sec:conclusions}.

\begin{small}
    \begin{tabular}{clcc}
    %\centering
    \multicolumn{3}{c}{\textbf{Nomenclature}}\\
        \hline
        Symbol & Explanation & & Unit \\
        \hline
        \multicolumn{2}{l}{\textbf{\emph{Mathematical notation}}}\\
        $\mathcal{T}=\{0,\ldots, T-1\}$ & Set of time periods $t$ in a given time horizon  & & $-$ \\
        $\tau_s$ & Duration of a time period (sampling time) & & $\textnormal{[h]}$ \\
        $\I(\underline t,\overline t)=[\underline t,\overline t)\subseteq \mathcal{T}$ & Generic time interval  & & $-$ \\
                \multicolumn{2}{l}{\textbf{\emph{Model variables}}}\\
     %   $U$ & Number of entities of the community & & $-$ \\
        $\mathcal{U}$ & Set of producer entities with storage, $\mathcal{U}=\{1,\ldots, U\}$ & & $-$ \\
     %   $T$ & Number of time slots of the optimization horizon & & $-$ \\ 
     %   $R$ & Number of DR requests in a day & & $-$ \\
        %$\underline{S}_u$ & Minimum capacity of energy storage of entity $u$ & & $\textnormal{[kWh]}$ \\
        $\overline{S}_u$ & Maximum BESS capacity of entity $u$ & & $\textnormal{[kWh]}$ \\
        $\overline{E}_u^{d}$ & Maximum BESS discharging energy of entity $u$ per time slot& & $\textnormal{[kWh]}$ \\
        $\overline{E}_u^{c}$ & Maximum BESS charging energy of entity $u$ per time slot& & $\textnormal{[kWh]}$ \\
        $\eta_u^{d}$ & Discharging efficiency of BESS of entity $u$ & & $-$ \\
        $\eta_u^{c}$ & Charging efficiency of BESS of entity $u$ & & $-$ \\
        $S_u^{0}$ & Energy level of BESS of entity $u$ at time  $t = 0$ & & $\textnormal{[kWh]}$ \\
        $S_u^{T}$ & Energy level of BESS of entity $u$ at time  $t = T$ & & $\textnormal{[kWh]}$ \\
        %$\overline{E}_u^{g}$ & Maximum energy injected into the grid by entity $u$ & & $\textnormal{[kWh]}$ \\
        $\pi_u^{g}(t)$ & Unitary price of energy sold to the grid by entity $u$ & & $\textnormal{[\euro/kWh]}$ \\
        $\pi_u^{s}$ & Unitary cost for operating the BESS of entity $u$ & & $\textnormal{[\euro/kWh]}$ \\
        $E_u(t)$ & Energy production of entity $u$ & & $\textnormal{[kWh]}$ \\
        $\hat E_u(t)$ & Forecast of $E_u(t)$ & & $\textnormal{[kWh]}$ \\
        $E_u^{g}(t)$ & Energy sold to the grid by entity $u$ & & $\textnormal{[kWh]}$ \\
        $E_u^{d}(t)$ & BESS discharging energy of entity $u$ & & $\textnormal{[kWh]}$ \\
        $E_u^{c}(t)$ & BESS charging energy of entity $u$ & & $\textnormal{[kWh]}$ \\
        $S_u(t)$ & Energy level of BESS of entity $u$ & & $\textnormal{[kWh]}$ \\
        $E^l(t)$ & Overall REC loads & & $\textnormal{[kWh]}$ \\
        $E^{n}(t)$ & Net energy injected into the grid by the REC & & $\textnormal{[kWh]}$ \\
        $E^{p}(t)$ & Overall energy generation by non-schedulable producers & & $\textnormal{[kWh]}$ \\
        $J_{u,0}$ & Net profit of entity $u$ when operating standalone  & & $\textnormal{[\euro]}$ \\
        $J_{u}$ & Net profit of entity $u$ when operating within the REC  & & $\textnormal{[\euro]}$ \\
        $\mathcal{R}_j$ & DR request & & $-$ \\
        $\I(\underline t_j,\overline t_j)$ & Time horizon of DR request ${\cal R}_j$ & &$-$\\
        $\underline{E}_j^{DR}, \overline{E}_j^{DR}$ & Lower and upper energy bounds for DR request ${\cal R}_j$ & &$\textnormal{[kWh]}$\\
        $\overline{\gamma}_{j}$ & Maximum DR reward for DR request $\mathcal{R}_j$ & & $\textnormal{[\euro]}$ \\
        $\gamma_{j}$ & DR reward for DR request $\mathcal{R}_j$ & & $\textnormal{[\euro]}$ \\
        $E_j^{DR}$  & Net energy injected into the grid in  $\I(\underline t_j,\overline t_j)$ & & $\textnormal{[kWh]}$ \\
        $\mathcal{R}$ & DR program, $\mathcal{R}=\{\mathcal{R}_1,\ldots, \mathcal{R}_R\}$ & & $-$ \\
        $\gamma$ & Total reward for DR program $\mathcal{R}$ & & $\textnormal{[\euro]}$ \\
        $\alpha$ & Fraction of $\gamma_j$ retained by the REC manager& & $-$ \\
        $\Psi_u$  & Revenue of entity $u$ for selling energy to the grid & & $\textnormal{[\euro]}$ \\        
        $\xi_u$  & DR reward assigned to entity $u$ & & $\textnormal{[\euro]}$ \\
        $z_{r,1}, z_{r,2}, z_{r,3}$ & Binary variables & & $\textnormal{[-]}$ \\
        \hline
    \end{tabular}
\end{small}

\section{Problem formulation and modeling}\label{sec:modeling}
The REC considered in this paper consists of a set of entities (participants) $\cU$, each equipped with a renewable generator, e.g., a PV plant, and a battery energy storage system, whose operation is scheduled by the REC manager via a centralized energy controller. The REC is also assumed to include additional entities composed by pure generators not connected to a BESS, as well as entities represented by  loads. Such players contribute to the REC energy balance as a whole, but are not subject to scheduling. 

\subsection{REC entity model}\label{ss:Em}
Operation decisions are assumed to be taken at discrete time instants $t$ within a given time horizon ${\cal T}=\{0,\dots,T-1\}$, e.g., one day. For each $t\in\cal T$, let $E_u(t)$ represent the amount of energy generated by entity $u\in\cU$ in the corresponding time slot, i.e., in the time frame beginning at $t$ and ending at $t+1$. Moreover, let $E^c_u(t)$ and $E^d_u(t)$ denote the controlled variables representing the energy injected into and drawn from the BESS, respectively, during time slot $t$, and let $S_u(t)$ be the BESS energy level at the beginning of the time slot, whose dynamics is modeled by the difference equation
\begin{equation}
	S_u(t+1) = S_u(t) +\eta_u^{c}E_u^{c}(t) - \frac{1}{\eta_u^{d}}E_u^{d}(t),
\end{equation}
where $0<\eta^c_u<1$ [$0<\eta^d_u<1$] represents the charging [discharging] efficiency. The controlled variables $E_u^{c}(t)$ and $E_u^{d}(t)$ are assumed to be bounded, i.e.,
\begin{equation}
	0 \leq E_u^{c}(t) \leq \overline{E}_u^{c}, \quad
	0 \leq E_u^{d}(t) \leq \overline{E}_u^{d} ,
\end{equation}
while $S_u(t)$ is bounded by the storage capacity, i.e., 
\begin{equation}
	0  \leq  S_u(t)	 \leq \bar S_u.
\end{equation}
Let $E_u^{g}(t)$ be the amount of energy injected into the grid in time slot $t$. The energy balance of participant $u$ relative to such time slot is therefore expressed by
\begin{equation}
	E_u^{g}(t) = E_u(t)-E_u^{c}(t)+E_u^{d}(t).
\end{equation}
The energy amount $E_u^g(t)$ is sold according to a known pricing signal $\pi_u^g(t)$ assumed known in advance, while the BESS is subject to a unitary operation cost $\pi_u^s$. Then, the net profit obtained by entity $u$ from energy sale over the time horizon $\T$ is given by
\begin{equation}
	J_{u,0} =\sum_{t\in\T} \Big[\pi_u^{g}(t)E_u^{g}(t)-\pi_u^{s}\big(\eta_u^{c}E_u^{c}(t) +\frac{1}{\eta_u^{d}}E_u^{d}(t)\big)\Big] .
\end{equation}
Finally, the total energy provided to the REC by non-schedulable producers is denoted by $E^p(t)$, while the overall load is indicated with $E^l(t)$, so that the net energy injected into the grid by the REC in time slot $t$ reads
\begin{equation}\label{eq:toten}
    E^{n}(t) = \sum_{u\in\U}  E^g_u(t) \ + E^p(t) - E^l(t) .
\end{equation}

\subsection{Demand response model}\label{ss:DRm}
\begin{figure}[h]
	\centering
	\includegraphics[width=0.5\linewidth]{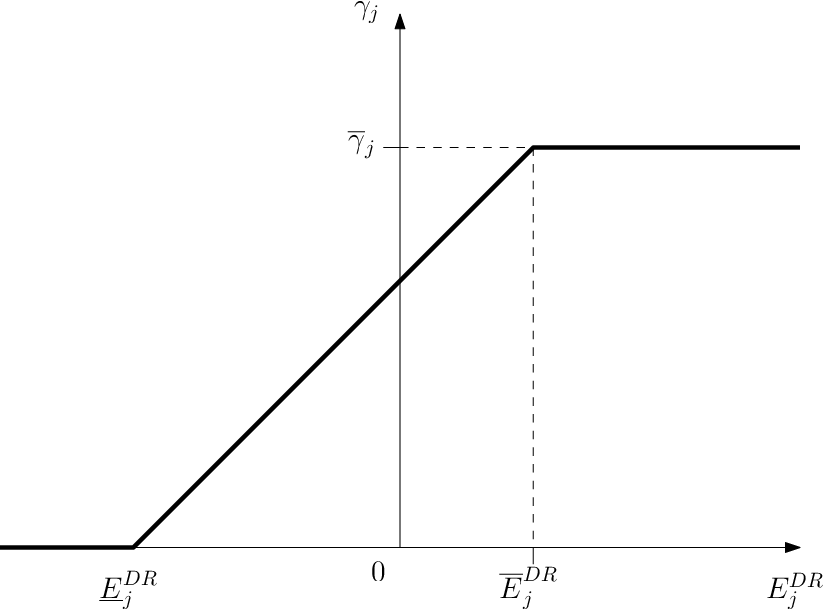}
	\caption{Overall reward related to the $j-$th DR request as a function of the net energy injected in the grid.}
	\label{fig:DR_DSO}
\end{figure}

The following demand response model based on price-volume signals is considered in this paper. A DR program $\cal R$ is modeled as a sequence of DR requests sent out by the DSO to the REC manager within the time frame $\cal T$, each one consisting of a time horizon and an associated monetary reward function. A suitable reward is granted to the REC if the net energy injected into the grid by the REC falls within suitable bounds. More specifically, a DR request ${\mathcal R}_j$ is defined by the following tuple:
\begin{equation}
\mathcal{R}_j=\left\{\I(\underline t_j,\overline t_j),\ \underline{E}_j^{DR},\  \overline{E}_j^{DR},\ \overline{\gamma}_j \right\} ,
\end{equation}
where $\I(\underline t_j,\overline t_j)\subseteq \cal T$ is the time horizon and $\underline{E}_j^{DR}, \overline{E}_j^{DR}$, and $\overline{\gamma}_j$ are positive bounds.
Let
\begin{equation}\label{eq:edr}
    E_j^{DR} = \sum_{t\in\I(\underline t_j,\overline t_j)}  E^{n}(t)  
\end{equation}
be the net energy injected into the grid by the REC within the time frame $\I(\underline t_j,\overline t_j)$. Then, the reward $\gamma_j$  corresponding to $\mathcal{R}_j$ granted to the REC is given by (see Fig.~\ref{fig:DR_DSO})
{\begin{equation}\label{eq:reward}
\gamma_j=
\begin{cases}
	\overline{\gamma}_j & \text{if } E_j^{DR} \geq \overline{E}_j^{DR}\\[2mm]
	\frac{(E_j^{DR} - \underline{E}_j^{DR})}{(\overline{E}_j^{DR} - \underline{E}_j^{DR})}\overline{\gamma}_j & \text{if } \underline{E}_j^{DR} \leq E_j^{DR} \leq \overline{E}_j^{DR} \\[2mm]
	0 & \text{if } E_j^{DR} \leq \underline{E}_j^{DR}
\end{cases} .
\end{equation}}

The overall reward associated to a given DR program $\cal R$ is given by
\[
\gamma = \sum_{j:\mathcal{R}_j\in\cR} \gamma_j .
\]
It is assumed that the REC manager redistributes a portion 
$$\xi = \alpha \gamma ~~~(0<\alpha<1)$$
of the overall reward $\gamma$ among the entities in $\cal U$ according to a fairness policy introduced later on, while retaining the remaining fraction to cover the profit of the REC manager itself and the reward granted to other REC participants (i.e., loads and generators without storage).

\section{Optimal REC management under DR program}\label{sec:optimization}
In this section we present the main contribution of this paper, which consists in the design of a three-step optimization procedure aimed at optimally operating the BESS resources of the community in the presence of DR programs. 
The key feature of the proposed method is to guarantee that for each $u\in\cal U$, the total profit obtained by joining the REC is always  greater or equal to the maximum achievable profit from energy sales obtained by optimally managing the BESS resources in an autonomous fashion, which is derived as a baseline in the first step. In the second step, an optimal scheduling of the energy storage systems of producers is computed so that a suitable performance index is maximized, accounting for the constraints arising from the models sketched in the previous section. In the last step, a policy for the distribution of the DR rewards among the participants is devised based on a fairness principle.

Given a scheduling time horizon $\mathcal{T}$, the proposed three-step design procedure is broken down as follows.
\begin{itemize}[leftmargin=1.25cm]
\item[\sl Step 1.] The optimal perspective profit $J_{u,0}^{*}$ that each single entity $u\in\U$ can achieve over $\mathcal{T}$ without joining the REC, i.e., without participating in the DR program, is calculated. The overall profit  $J_{0}^{*}$ of all $u\in\cal U$ is computed and exploited in the second step as a lower bound constraint on the overall expected profit when participating in the REC.
\item[~~\sl Step 2.] The optimal scheduling of the control variables of the BESSs of the whole REC under the DR program is computed in order to maximize a suitable overall REC performance index $H$ while guaranteeing an overall profit of at least $J_{0}^{*}$. Possible choices for $H$ are discussed at the end of this section.
\item[~~\sl Step 3.] The redistribution of the DR reward obtained using the control policy in Step 2 is computed, ensuring both fairness among participants and a total (i.e., energy sales plus DR rewards) revenue $J^*_u\geq J^*_{u,0}$ for each $u\in\cal U$.
\end{itemize}

Let us define
\begin{equation}\label{eq:defpsi}
   \Psi_u = \sum_{t\in\T} \Big[\pi_u^{g}(t)E_u^{g}(t)-\pi_u^{s}\big(\eta_u^{c}E_u^{c}(t) +\frac{1}{\eta_u^{d}}E_u^{d}(t)\big)\Big] ,
\end{equation}
   which represents the net operation profit of entity $u$ over the time horizon $\T$ arising from energy sales. Moreover, let $\hat E_u(t)$, $t\in\T$ denote a forecast of the renewable energy generation of $u$.
Step~1 of the above procedure can be accomplished by solving for each $u\in\cal U$ an optimization problem involving the set of decision variables
\begin{equation}\label{eq:theta}
    \Theta_{u} = \left\{\{E_u^g(t),\ E_u^c(t),\ E_u^d(t),\ S_u(t),\ \forall t\in\T\},\ S_u(T) \right\}
\end{equation}
and formulated as follows.
\begin{problem}\label{pb:alone}
\begin{equation*}
	J_{u,0}^{*} = \max_{\Theta_{u}} \Psi_u
\end{equation*}
$\qquad$~~subjected to:
\begin{eqnarray}
	0 \leq E_u^{c}(t) \leq \overline{E}_u^{c}, \quad
	0 \leq E_u^{d}(t) \leq \overline{E}_u^{d} & ~ \label{eq:constr_begin} \\
	S_u(t+1) = S_u(t) +\eta_u^{c}E_u^{c}(t) - \frac{1}{\eta_u^{d}}E_u^{d}(t)  & ~\\
	0 \leq  S_u(t)	 \leq \overline{S}_u & ~~~~~~~~~\forall t\in\T \\
	%\quad 0  \leq {E_u^{g}(t)} \leq \overline{E}_u^{g} & ~~~~~~~~~\forall t\in\T \\
	E_u^{g}(t)+E_u^{c}(t)	= \hat E_u(t)+E_u^{d}(t) & ~ \label{eq:constr_extra1} \\
    E_u^{c}(t) E_u^{d}(t)	= 0 & ~ \label{eq:compl_constr} \\
	E_u^{c}(t) \leq \hat E_u(t) & ~ \label{eq:constr_nogrid}
\end{eqnarray}
\begin{equation}\label{eq:constr_end}\vspace*{-0.2mm}
	S_u(0) =S_u^{0},\quad S_u(T) =S_u^{T}
\end{equation}
\end{problem}
In Problem \ref{pb:alone}, constraints \eqref{eq:constr_begin}-\eqref{eq:constr_extra1} derive from the models in Section \ref{ss:Em}, where the generation forecast time series $\hat E_u(t)$ is used in place of the actual generation $E_u(t)$.
Moreover, \eqref{eq:compl_constr} and \eqref{eq:constr_nogrid} are used to avoid simultaneous BESS charging/discharging and BESS charging from the grid, respectively. Finally, constraints \eqref{eq:constr_end} ensure that the initial energy level of the storage is reset to a prescribed value at the end of the operation horizon.
\begin{remark}
 It is not difficult to show that the complementarity constraint \eqref{eq:compl_constr} in Problem \ref{pb:alone} is always satisfied at the optimum and therefore it can be omitted. In fact, among all feasible solutions for which the term $\eta_u^{c}E_u^{c}(t) - \frac{1}{\eta_u^{d}}E_u^{d}(t)$ is constant, the ones satisfying either $E_u^{c}(t)=0$ or $E_u^{d}(t)=0$ yield a higher objective value. By neglecting constraint \eqref{eq:compl_constr}, Problem 1 turns out to be a linear program.
\end{remark}

In order to devise the optimal REC scheduling strategy in Step 2, it is first convenient to reformulate the reward policy \eqref{eq:reward} associated to each DR request $\mathcal{R}_j$ as a set of linear inequalities involving binary variables. Indeed, it is easily seen that $E_j^{DR}$ and $\gamma_j$ satisfy \eqref{eq:reward} if and only if there exist $z_{j,1},z_{j,2},z_{j,3}$ such that
\begin{eqnarray}
        z_{j,1},z_{j,2},z_{j,3} \in\{0,1\} \label{eq:binaryv} \\
	z_{j,1}+z_{j,2}+z_{j,3} \leq 1 \label{eq:binary_vars}\\ 
-M z_{j,1} + \underline{E}_j^{DR} z_{j,2} + \overline{E}_j^{DR} z_{j,3} \leq E_j^{DR} \leq \underline{E}_j^{DR} z_{j,1} + \overline{E}_j^{DR} z_{j,2} + M z_{j,3} \\
    -M(1 - z_{j,3}) \leq \gamma_j - \overline{\gamma}_j^{DR} \leq M(1 - z_{j,3}) \\
    -M(1 - z_{j,2})  \leq \gamma_j - \frac{(E_j^{DR} - \underline{E}_j^{DR})}{(\overline{E}_j^{DR} - \underline{E}_j^{DR})} \overline{\gamma}_j^{DR} \leq M(1 - z_{j,2}) \\
    -M(1 - z_{j,1})  \leq \gamma_j \leq M(1 - z_{j,1}) \label{eq:jlast}
\end{eqnarray}
where $M\gg0$ denotes a constant big enough to avoid inconsistencies in the formulation.

The proposed optimal storage scheduling strategy in Step 2 is computed via the solution of the following optimization problem, where
\[
J_0^* = \sum_{u\in\cU} J_{u,0}^{*} 
\]
and the set of decision variables is defined as
\begin{align}
\Theta = \left\{\Theta_{u}, \gamma_j,z_{j,1},z_{j,2},z_{j,3}, \forall u \in \cU,~ \forall j:{\cal R}_j \in \cR \right\}. 
\end{align}

\begin{problem}\label{pb:manager_new}
\begin{equation}
H^* =\max_{\Theta} ~H
\end{equation}
$\qquad$~~subjected to:
\begin{eqnarray}
\eqref{eq:defpsi},\eqref{eq:constr_begin}-\eqref{eq:constr_end},~ \forall u\in\cU,~ \forall t\in\T \\
\eqref{eq:toten},\eqref{eq:edr}, \eqref{eq:binaryv}-\eqref{eq:jlast},~ \forall j:{\cal R}_j\in\cR  \\
%\gamma^{DR}=\ds\sum_{j:\mathcal{R}_j\in\cR}\gamma_j^{DR} \\
%\Psi =\sum_{u\in\U} \Psi_u
\sum_{u\in\U} \Psi_u + \alpha \sum_{j:\mathcal{R}_j\in\cR} \gamma_j - J_{0}^{*}  \ge0 \label{eq:delta_ge0_new}
\end{eqnarray}
\end{problem}
\begin{remark}
To avoid complicating the notation, constraint \eqref{eq:toten} in Problem \ref{pb:manager_new} is assumed to be evaluated for $E^p(t)$ and $E^l(t)$ equal to suitable forecasts of the respective variables.
\end{remark}
Once Problem \ref{pb:manager_new} is solved, the optimal $\Theta_u^*$, $u\in\cU$ yields the optimal scheduling strategy of the BESS resources. Moreover, let $\gamma^*_j$ and $\Psi^*_u$ be the optimal values of $\gamma_j$ and $\Psi_u$, respectively. Notice that $\Psi_u^*$ represents the operation profit of entity $u$ under the optimal policy $\Theta^*_u$, not considering any extra reward from DR. 
Finally, the optimal total DR reward for all entities in $\cal U$ is given by
\begin{equation}\label{eq:DR_alpha}
\xi^*= \alpha \sum_{j:\mathcal{R}_j\in\cR} \gamma_j^* .
\end{equation}
Such a reward is to be shared among entities in $\cal U$ according to the policy described in the next subsection.
\begin{remark}
From a computational viewpoint, Problem \ref{pb:manager_new} turns out to be a mixed-integer linear program involving a number of binary variables equal to $3R$ being $R$ the cardinality of $\cal R$. Since the daily number of DR requests amounts to a few units at most, the total number of binary variables is low regardless of the community size. Notice that the observation concerning the possibility of neglecting the nonlinear constraint \eqref{eq:compl_constr} still stands for Problem \ref{pb:manager_new} if the objective function $H$ is reasonably defined, e.g., as a positively weighted sum of the profits $\Psi_u$ and of the DR rewards $\gamma_j$, due to the monotonicity of $\gamma_j$ with respect to $E_j^{DR}$. This will indeed be the case with the indices proposed later on.
\end{remark}

\subsection{DR reward redistribution}
Let $\xi^*_u$, $u\in \cal U$ denote a partition of $\xi^*$ among the entities of $\cal U$, i.e.,
\begin{equation}\label{eq:partieta}
 \sum_{u\in\cal U} \xi^*_u = \xi^*, ~~~ \xi^*_u\geq 0 .
 \end{equation}
The total profit of entity $u\in\cal U$ in the presence of DR is then given by 
$$J^*_u = \Psi_u^* + \xi^*_u.$$
As previously discussed, the partitioning is deemed convenient if for all $u\in\cal U$, the total profit equals at least the baseline $J^*_{u,0}$, i.e., 
\begin{equation}\label{eq:conveniente}
J_u^*  \geq J^*_{u,0}~~~\forall u\in\cal U .
\end{equation}
Hence, we look for a partition $\xi^*_u$, $u\in\cal U$ satisfying \eqref{eq:partieta},\eqref{eq:conveniente}. Many possible schemes can be devised to this purpose. In the following, we propose a reward assignment based on a fairness principle, which guarantees to all entities a reward proportional to their baseline profit. To this purpose, we introduce the ratio
\begin{equation}\label{eq:defrho}
\rho =\frac{\Psi^* + \xi^*  - J^*_{0}}{J^*_{0}}
\end{equation} 
being
\[
\Psi^* = \sum_{u\in\U} \Psi_u^* .
\]
By virtue of constraint \eqref{eq:delta_ge0_new} in Problem \ref{pb:manager_new}, it turns out that $\rho\geq 0$. Let us define $\xi_u^*$ as
\begin{equation}\label{eq:partibuona}
\xi_u^* = (1+\rho) J^*_{u,0} - \Psi^*_u .
\end{equation}
This corresponds to a total entity profit
\[
J^*_u = (1+\rho) J^*_{u,0} .
\]
Hence, the partition  \eqref{eq:partibuona} trivially satisfies \eqref{eq:conveniente}, while \eqref{eq:partieta} simply follows by taking the sum of both sides of \eqref{eq:partibuona} over all $u\in{\cal U}$ and using \eqref{eq:defrho}. \\
Clearly, the additional profit gained by $u$ by joining the REC is given by
\[
\delta_u = \rho J^*_{u,0} .
\]
Note that fairness among entities is guaranteed by the proposed sharing policy since the extra profit obtained by each producer $\delta_u$ is proportional to the profit achieved if operating autonomously. In other words, $\rho$ denotes the ratio between the additional profit and the baseline, for each entity. 

It is worth remarking that reward distribution policies different from the one proposed here can be easily devised. In fact, the computations in steps 1 and 2 are independent of the particular distribution strategy implemented in step 3.

\subsection{Performance indices}\label{sec:Performanceindices}
Depending on the main aims of the community (e.g., maximizing the REC manager profit, maximize the overall revenue of entities, etc.), several objective functions $H$ can be defined for Problem \ref{pb:manager_new}. The following two choices are proposed and compared in this paper:
\begin{itemize}
\item \textit{REC manager interest:} the objective function is taken equal to the revenue of the REC manager, i.e.,
\begin{equation}
H = H^{M}=\ds(1-\alpha) \sum_{j:\mathcal{R}_j\in\cR} \gamma_j.
\end{equation}
Note that maximizing $H^M$ is equivalent to maximizing the overall DR reward gained by the community.
\item \textit{Overall entities' interest:} the considered objective function is the total profit of the entities $u\in\U$, i.e.,
\begin{equation}
H = H^{E}=\sum_{u\in\U} \Psi_u + \alpha \sum_{j:\mathcal{R}_j\in\cR} \gamma_j .
\end{equation}
\end{itemize}
Different choices of the performance index $H$ are indeed possible, as well as different definitions of the reward redistribution policy. Such extensions will be the topic of further research.
%\begin{remark}
%It is worth noting that the inclusion of constraint \eqref{eq:delta_ge0_new} might seem unnecessary when the %cost function $J$ is set to $J^{E}$. However, we retain this constraint to ensure a unified formulation of %the optimization problem across different cost functions. Notably, when optimizing $J^E$ or $J^M$, situations %may arise where $\Psi_{\mathcal{U}}^{REC}$ falls below $J_{\mathcal{U}}^{0}$. In such cases, the shortfall is %effectively compensated by the rewards provided through the DR program.
%\end{remark}

\section{Test cases} \label{sec:test_cases}
This section offers two examples to validate the proposed approach. The first one is a simple illustrative example aimed at showing the main features of the procedure, while the second one has the purpose to evaluate the performance and the computational feasibility of the proposed technique in a larger scale practical setting. In all simulations, the optimization horizon $\cal T$ is assumed to span $24$ hours with a sampling time of $\tau_s=15$ minutes. 
In both examples, the parameters of the energy storage systems as well as the energy generation profiles of the various entities have been obtained by a suitable scaling of the real data used in \cite{stentati2023optimization}.
For all entities $u$, the considered initial and final BESS energy levels are $S_u^0=S_u^T=0$, the charging and discharging efficiencies are set to $\eta_u^{c}=\eta_u^{d}= 0.95$, and the unit price of the energy sold to the grid $\pi_u^{g}(t)$ is set as the time series depicted in green in Fig.~\ref{fig:ex1_alone}, while the unitary cost for operating the storage system is $\pi_u^{s}=0.01\textnormal{\euro/kWh}$. Simulations are performed assuming the availability of suitable forecasts of the generation profiles of each PV producer.
The generation data cover the period from April 1, 2019 to April 30, 2019, for a total of $30$ days. Simulations and optimization of Problems \ref{pb:alone} and \ref{pb:manager_new} have been implemented in Python and solved using the CPLEX solver \cite{cplex} on an Intel i7-11700@3.60 GHz, 16 GB of RAM. Results from the simulations are discussed in Section \ref{sec:discussion}.

\begin{example} 
{\rm In this illustrative example, a community composed of two PV producers equipped with a BESS is considered. 
%works and to better illustrate the obtained results.
The nominal peak power of the PV plant of the first entity is set to $\overline P_1=700\textnormal{kW}$, with a battery capacity $\overline S_1=500\textnormal{kWh}$. In contrast, the second entity can supply generation with
$\overline P_2=300\textnormal{kW}$ peak power, with a battery capacity 
$\overline S_2=250\textnormal{kWh}.$
We assume two DR requests for each day. 
%The maximum reward $\overline{\gamma}_j$, time intervals $\underline{t}_j$ and $\overline{t}_j$, and energy limits $\overline{E}_j^{DR}$, $\underline{E}_j^{DR}$ 
The parameters associated with each request are provided in Table \ref{table:energy_requests}. Such DR requests are assumed to remain the same in all the considered days.
 
In Table \ref{tab:table2}, the following daily amounts are reported for 5 simulation days: the optimal total daily profit of entities when operating outside the REC ($J^*_{0}$), the total daily extra profit of the entities when joining the REC ($\delta_1+\delta_2$), and the total daily DR reward received at REC level ($\gamma_1^{*}+\gamma_2^{*}$). Results for both objective functions $H^{E}$ and $H^{M}$ introduced in Section \ref{sec:Performanceindices} are reported.
The individual profit for each entity under Problem \ref{pb:alone} ($J^*_{u,0}$) and Problem \ref{pb:manager_new} ({$J^*_{u}$}) for both objective functions, as well as the corresponding $\rho$ are reported in Table \ref{tab:table3}. The time evolution of relevant energy signals concerning Entity 1 (along with that of energy selling price) in a representative day (day 24) are reported in Fig. \ref{fig:ex1_alone}, assuming the entity operates individually without joining the REC (Problem \ref{pb:alone}). The case of the same entity participating in the REC according to the proposed scheduling under the objective function $H^E$ is depicted in Fig. \ref{fig:ex1_REC}.
The daily extra profits gained by both entities under the objective functions $H^E$ and $H^M$, are illustrated in Fig. \ref{fig:ex1_delta}. 
Fig. \ref{fig:ex1_DR} shows the daily reward for the two DR requests $\mathcal{R}_1$ and $\mathcal{R}_2$, under the two objective functions $H^E$ and $H^M$. 

\begin{table}[!t]
\centering
\begin{tabular}{|c|c|c|c|c|c|c|}
\hline
Request ID & $\underline t_j$ & $\overline t_j$&$\underline{E}_j^{DR}\textnormal{[kWh]}$&$\overline{E}_j^{DR}\textnormal{[kWh]}$ & $\overline{\gamma}_j\textnormal{[\euro]}$&$\alpha$ \\
\hline
${\cal R}_1$          & 08:00      & 09:00 &0  & 800          & 65&0.85     \\
${\cal R}_2$         & 17:00      & 18:00   &0 & 1400         & 65&0.85   \\
\hline
\end{tabular}
\caption{Example 1. Set of DR requests for all days.}
\label{table:energy_requests}
\end{table}

\begin{table*}[!t]
    \centering
    \centering
     \small % Reduce the font size of the table
    \renewcommand{\arraystretch}{1.2} % Adjusts the height of all rows
    \setlength{\tabcolsep}{2pt} % Adjust column separation
    \begin{tabular}{||c||c||c|c||c|c||c|c|c||}
        \cline{3-6}  
        \multicolumn{2}{c|}{}&
    \multicolumn{2}{c||}{$ H = H^{E}$} & \multicolumn{2}{c||}{ $H = H^{M} $} \\
        \cline{1-6}
           {Date} & {$J^*_{0}$[\euro]} &   $\delta_1+\delta_2$[\euro]&  $\gamma_1^{*}+\gamma_2^{*}$[\euro] &$\delta_1+\delta_2$[\euro]&  $\gamma_1^{*}+\gamma_2^{*}$[\euro] \\
           %{Date} & {$J^*_{0}$[\euro]} &   $\displaystyle \sum_{u \in \mathcal{U}}\delta_u$[\euro]&  $\ds\sum_{j:\mathcal{R}_j\in\cR}\gamma_j^{*}$[\euro] &$\displaystyle\sum_{u \in \mathcal{U}}\delta_u$[\euro]&  $\ds\sum_{j:\mathcal{R}_j\in\cR}\gamma_j^{*}$[\euro]  \\           
        \hline
        04-01 &  438.33  & 40.59 & 74.56 & 13.84& 82.15\\
        \hline
        04-02 & 355.10  & 23.63 & 47.48 & 2.08 &62.49\\
        \hline
        04-03 & 403.71   & 34.39 & 67.28 & 12.35 & 77.46 \\
        \hline
        04-04 &  406.58 & 32.68 & 64.60 & 15.65& 75.59\\
        \hline
        04-05 &  375.59 & 29.11 & 61.20 & 0.58 & 72.08\\
        \hline
    \end{tabular}
    \caption{Example 1. Total daily amounts related to two different objective functions $H$, in 5 simulated days.}
    \label{tab:table2}
\end{table*}

\begin{table*}[!t]
    \centering
     \small % Reduce the font size of the table
    \renewcommand{\arraystretch}{1.2} % Adjusts the height of all rows
    \setlength{\tabcolsep}{2pt} % Adjust column separation
    \begin{tabular}{||c||c||c||c|c|c||c|c|c||}
           \cline{4-9}  
        \multicolumn{3}{c|}{}& 
        \multicolumn{3}{c||}{$H = H^{E}$} & \multicolumn{3}{c||}{ $ H = H^{M}$} \\
        \cline{1-9}
       {Date} & {$J^*_{1,0}$[\euro]} & {$J^*_{2,0}$[\euro]} & $J^*_{1}$[\euro]& $J^*_{2}$[\euro] & $\rho$&$J^*_{1}$[\euro]&  $J^*_{2}$[\euro]   &  $\rho$  \\
        \hline
        04-01 &306.48 & 131.85 & 334.86 & 144.06&0.092& 316.15& 136.01& 0.0316\\
        \hline
        04-02 &248.55 &106.55  & 265.09 & 113.64 &0.066& 250.01 &107.18&0.0050\\
        \hline
        %04-03 & 282.25 &121.46  & 306.29 & 131.81 &0.085& 290.88 & 125.18&0.00306\\
        04-03 & 282.25 &121.46  & 306.29 & 131.81 &0.085& 290.88 & 125.18&0.0031\\
        \hline
        04-04 &  284.28&122.30  & 307.13 &	132.13 &0.080& 295.233&127.01&0.0385 \\
        \hline
        04-05 &  262.87& 112.72 & 283.24& 121.46&0.077& 263.27&	 112.9&0.0015\\
        \hline
    \end{tabular}
    \caption{Example 1. Profit of each entity working outside the REC and profit and fairness index for different objective functions $H$, in 5 simulated days.}
    \label{tab:table3}
\end{table*}

\begin{figure*}[!t]
 \centering
 \includegraphics[width=1\textwidth, height=.3\textheight]{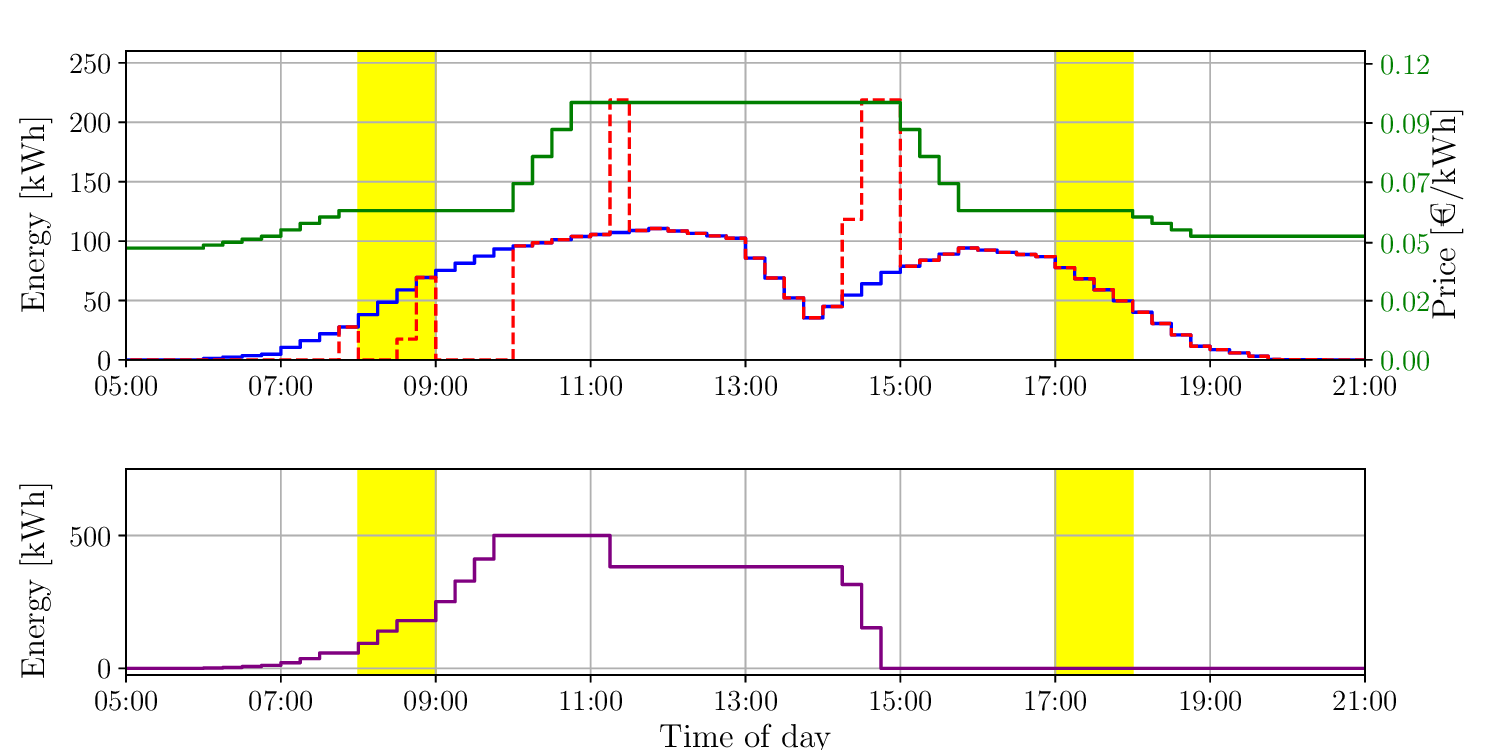}    
 \caption{Example 1. Results for Entity 1 on day 24 under Problem \ref{pb:alone}. Energy production forecast $\hat E_u(t)$ (blue), sold energy $E^g_u(t)$ (dashed red), energy selling price $\pi^g_u(t)$ (green), storage energy level $S_u(t)$ (purple), and time periods of DR requests (yellow).}
 \label{fig:ex1_alone} 
\end{figure*}

\begin{figure*}[!t]
 \centering
 \includegraphics[width=1\textwidth, height=.3\textheight]{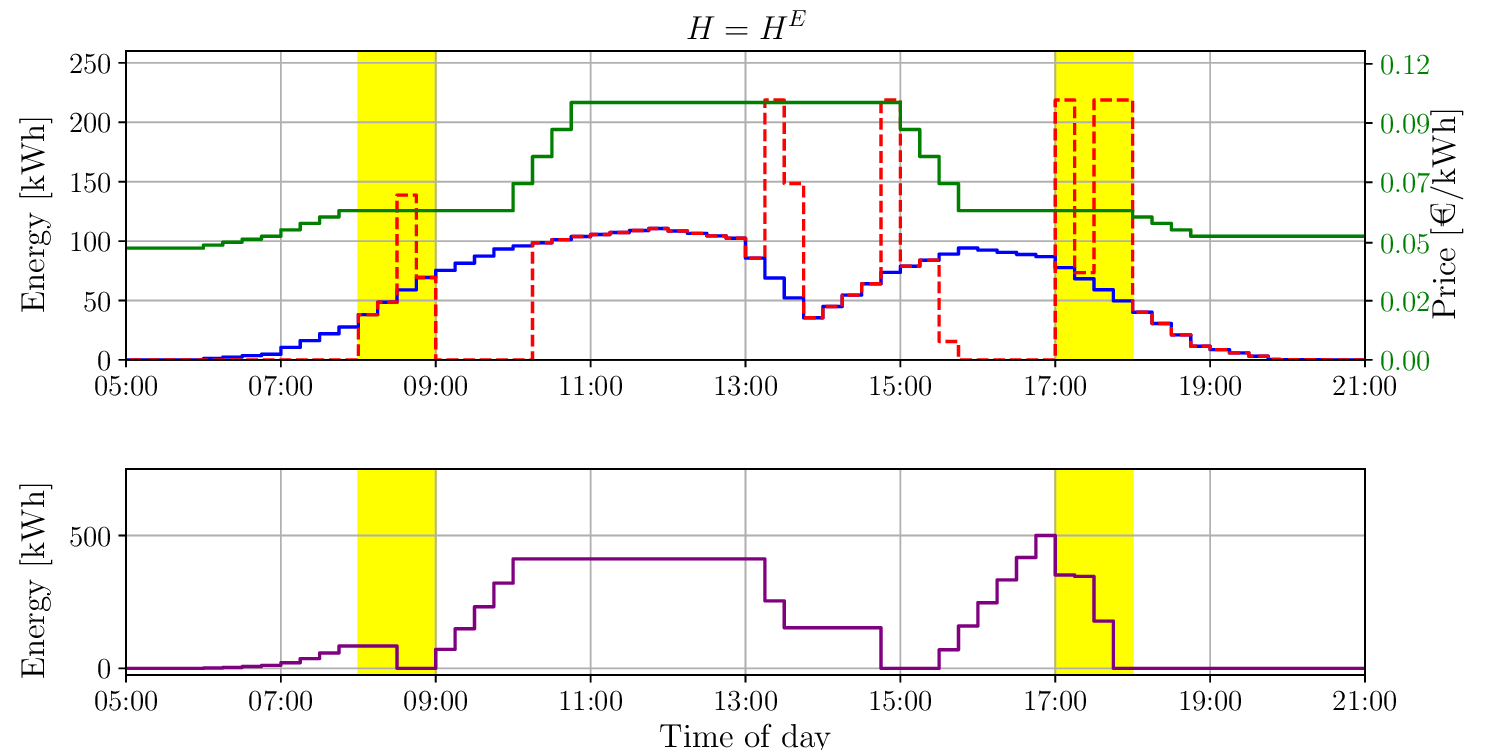}    
 \caption{Example 1. Results for Entity 1 on day 24 under Problem 2 with $H=H^{E}$. Energy production forecast $\hat E_u(t)$ (blue), sold energy $E^g_u(t)$ (dashed red), energy selling price $\pi^g_u(t)$ (green), storage energy level $S_u(t)$ (purple), and time periods of DR requests (yellow).}
 \label{fig:ex1_REC} 
\end{figure*}

\begin{figure*}[!t]
    \centering
    \includegraphics[width=\textwidth, height=.16\textheight]{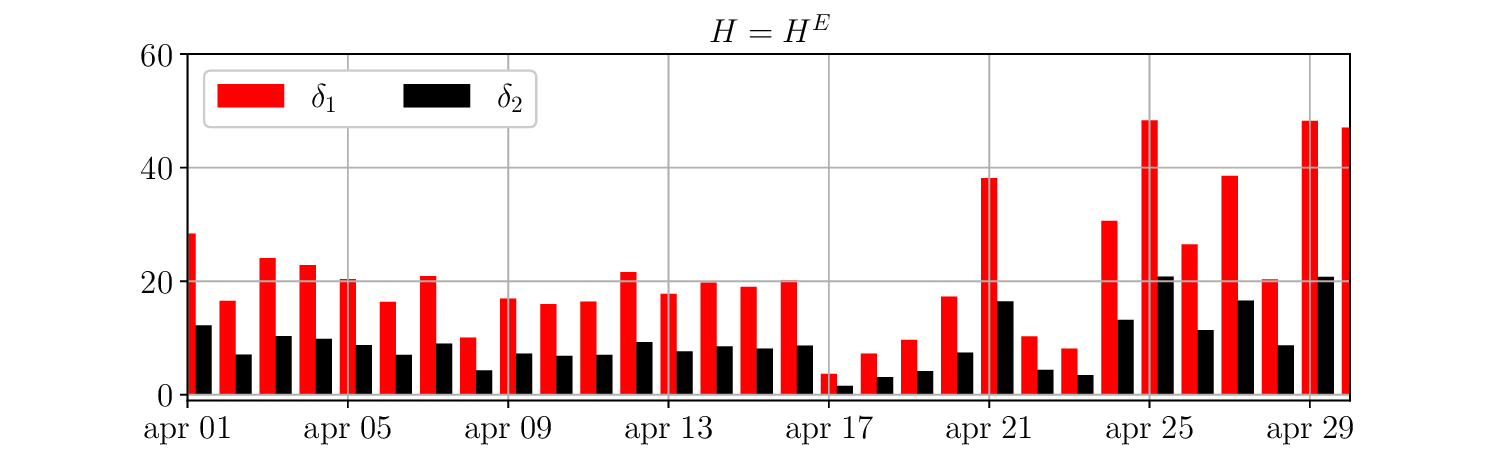}\\[1mm]
    \includegraphics[width=\textwidth, height=.16\textheight]{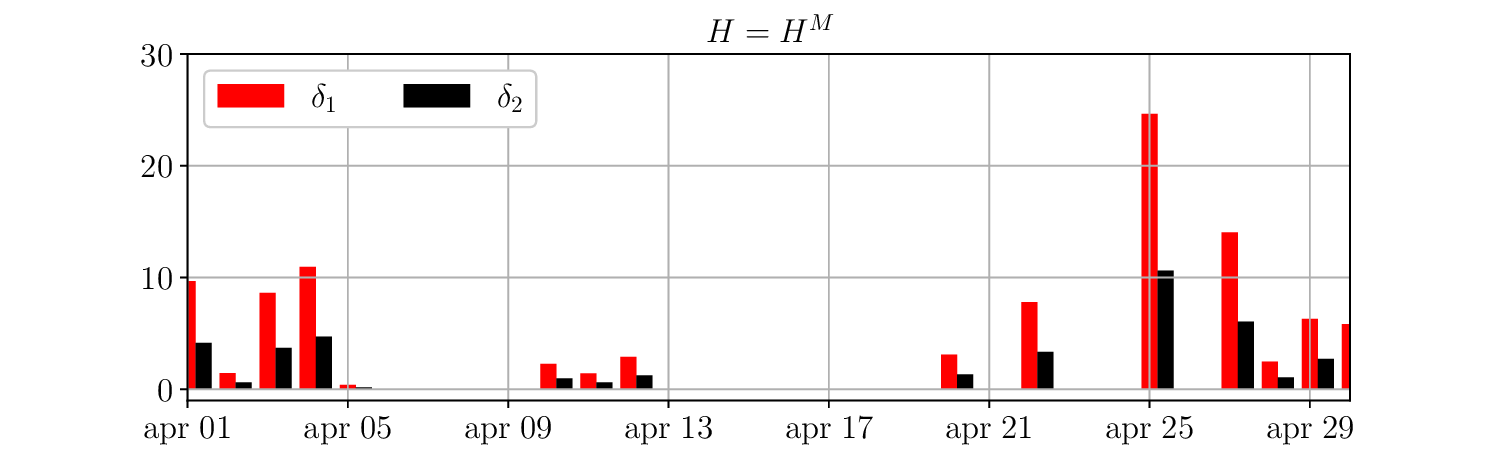}
    \caption{Example 1. Daily additional profits $\delta_1$ and $\delta_2$ under the objective functions $H^E$ (top) and $H^M$ (bottom).}
    \label{fig:ex1_delta}
\end{figure*}

\begin{figure*}[!t]
    \centering
    \includegraphics[width=\textwidth, height=.16\textheight]{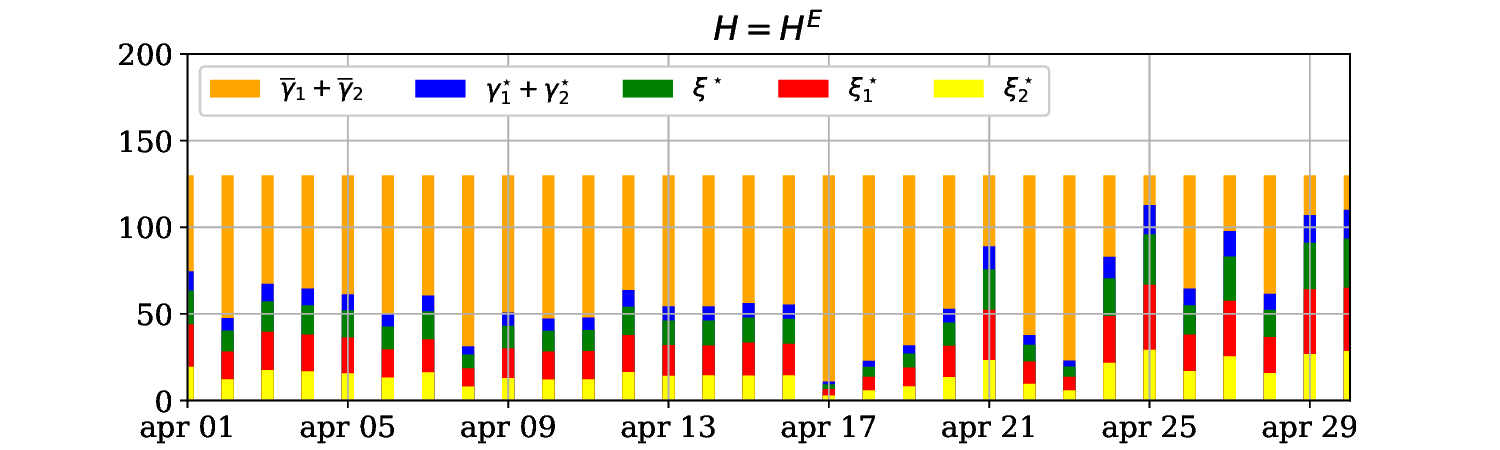}\\[1mm]
    \includegraphics[width=\textwidth, height=.16\textheight]{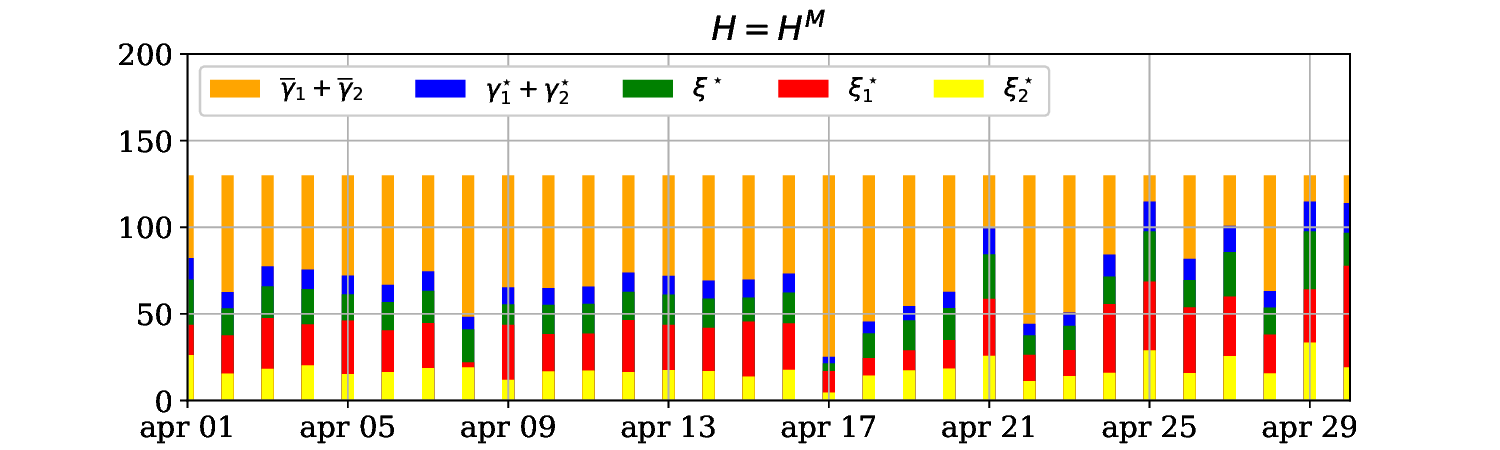}
    \caption{Example 1. Daily DR reward profiles for requests 1 and 2, under the objective functions $H^E$ (top) and $H^M$ (bottom). Maximum achievable DR reward $\overline{\gamma}_1+\overline{\gamma}_2$ (orange), actual DR reward at community level $\gamma_1^*+\gamma_2^*$ (blue), DR reward of Entity 1 $\xi_1^*$ (red) and 2 $\xi_2^*$ (yellow), and total reward of entities $\xi^*=\xi_1^*+\xi_2^*$ (green).}
    \label{fig:ex1_DR}
\end{figure*}
}

\begin{table*}[!t]
    \centering
    \small % Reduce the font size of the table
    \renewcommand{\arraystretch}{1.2} % Adjusts the height of all rows
    \setlength{\tabcolsep}{2pt} % Adjust column separation
    \begin{tabular}{|c|c|c|c||c|c||c|c||c|c||}
        \cline{3-10}  
        \multicolumn{2}{c|}{} & \multicolumn{4}{c||}{$\overline{\gamma}_j =40$[\euro]} & \multicolumn{4}{c||}{$\overline{\gamma}_j=100$[\euro]} \\
        \cline{3-10}
        \multicolumn{2}{c|}{} & \multicolumn{2}{c||}{$H = H^{E}$} & \multicolumn{2}{c||}{$H = H^{M}$} & \multicolumn{2}{c||}{$H = H^{E}$} & \multicolumn{2}{c||}{$H = H^{M}$} \\
        \cline{1-10}
        {Date} & {$J^*_{0}$[\euro]} &  $\delta_1+\delta_2$[\euro]&  $\gamma_1^{\star}+\gamma_2^{\star}$[\euro] & $\delta_1+\delta_2$[\euro]&  $\gamma_1^{\star}+\gamma_2^{\star}$[\euro] & $\delta_1+\delta_2$[\euro]&  $\gamma_1^{\star}+\gamma_2^{\star}$[\euro] & $\delta_1+\delta_2$[\euro]&  $\gamma_1^{\star}+\gamma_2^{\star}$[\euro] \\        
%        {Date} & {$J^*_{0}$} &  $\displaystyle \sum_{u \in \mathcal{U}}\delta_u$&  $\ds\sum_{j:\mathcal{R}_j\in\cR}\gamma_j^{\star}$ & $\displaystyle\sum_{u \in \mathcal{U}}\delta_u$ & $\ds\sum_{j:\mathcal{R}_j\in\cR}\gamma_j^{\star}$ & $\displaystyle \sum_{u \in \mathcal{U}}\delta_u$&  $\ds\sum_{j:\mathcal{R}_j\in\cR}\gamma_j^{\star}$ & $\displaystyle\sum_{u \in \mathcal{U}}\delta_u$ & $\ds\sum_{j:\mathcal{R}_j\in\cR}\gamma_j^{\star}$ \\
        \hline
        04-01 & 438.33 & 19.13 & 30.2 & 0.0 & 50.55 & 77.77 & 126.38 & 57.16 & 126.38 \\
        \hline
        04-02 & 355.10 & 10.08 & 18.1 & 0.0 & 38.45 & 51.41 & 96.13 & 37.37 & 96.13 \\
        \hline
        04-03 & 403.71 & 15.01 & 27.31 & 0.0 & 47.67 & 69.29 & 119.17 & 44.21 & 119.17 \\
        \hline
        04-04 & 406.58 & 14.13 & 26.16 & 3.59 & 46.52 & 66.68 & 116.29 & 62.90 & 116.29 \\
        \hline
        04-05 & 375.59 & 11.68 & 24.00 & 0.55 & 44.36 & 61.50 & 110.89 & 56.74 & 110.89 \\
        \hline
    \end{tabular}
    \caption{Example 1. Extra profits and achieved DR rewards under $H^E$ and $H^M$, for different values of the maximum DR reward.}
    \label{tab:4}
\end{table*}
\end{example}

\begin{figure*}[!ht]
    \centering
    \includegraphics[width=\textwidth, height=.18\textheight]{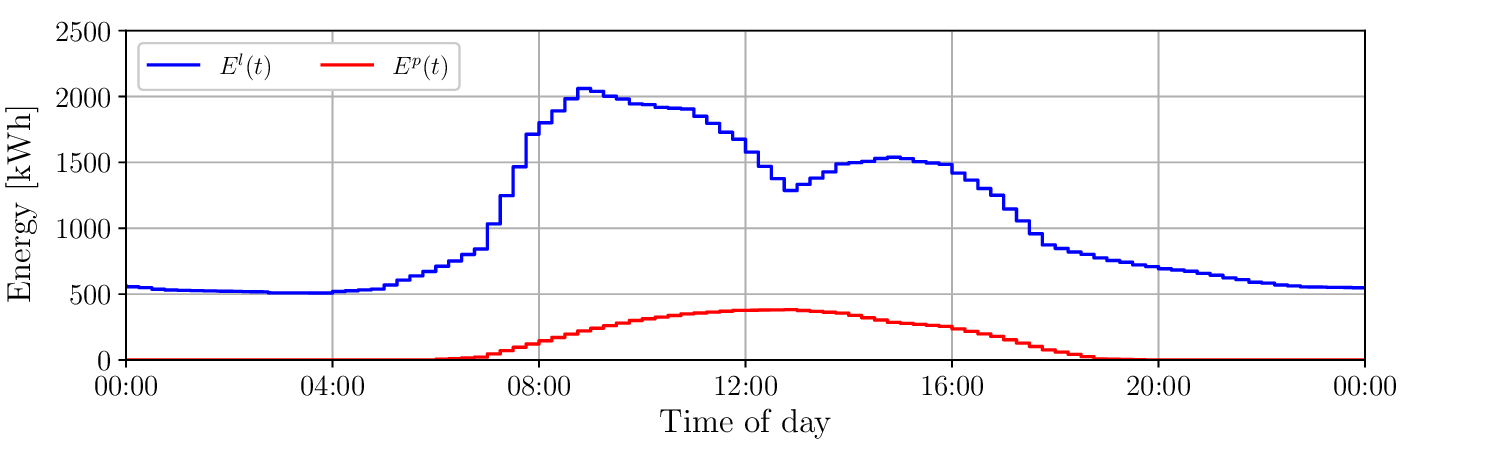}
    \caption{Example 2. Overall REC loads (blue), overall energy generation by non-schedulable producers (red) on day 3.}
    \label{fig:ex2_load}
\end{figure*}

\begin{figure*}[!ht]
    \centering
    \includegraphics[width=\textwidth, height=.18\textheight]{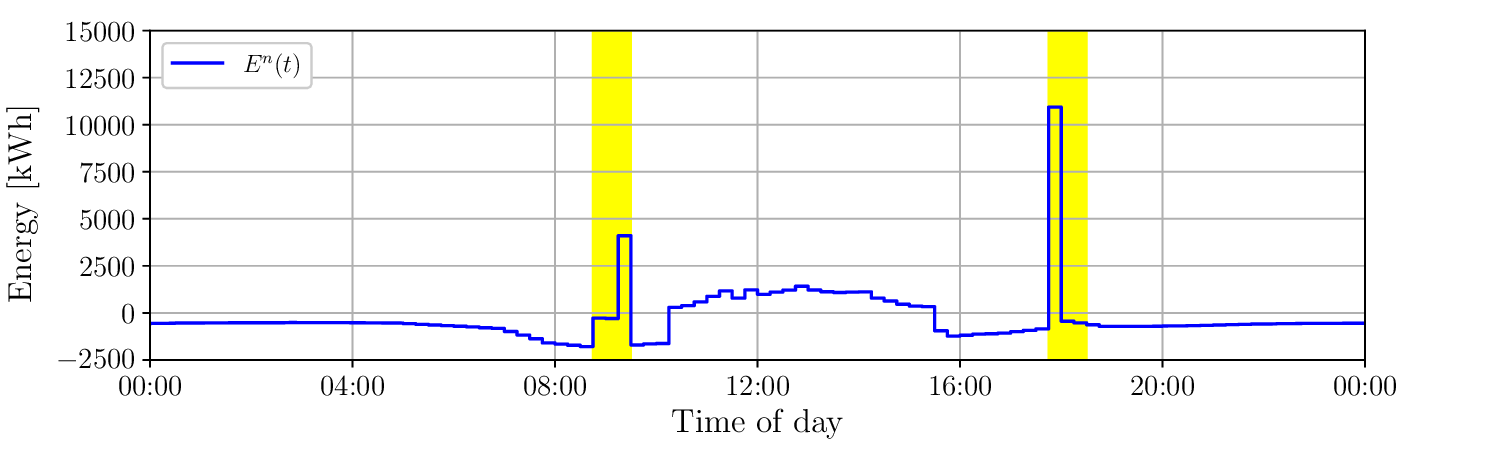}
    \caption{Example 2. Net energy injected into the grid by the REC  under $H^E$ (blue), and period of DR requests (yellow) on day 3.}
    \label{fig:ex2_NET}
\end{figure*}

\begin{example} 
{\rm
In this example, a community composed of 30 PV producers equipped with energy storage systems is considered. The peak power $\overline{P}_u$ of the PV plant of each entity and the corresponding battery capacity $\overline{S}_u$ are shown in Table \ref{tab:ex2_parameters}.
%$$
%\overline{S}_u=\alpha_u^s
%\underbrace{ \max_{d \in D} \left\{ \sum_{t \in \tau} E_{u,d}(t) 
%\right\}}_{E_u^{\max}}, 
%$$ 
%where $E_{u,d}(t)$ is the generated energy by entity $u$, on day $d$ at time $t$, while $\alpha_u^s$ denotes a scaling factor randomly chosen between 0.1 and 0.25. The values of the parameters related to entities are detailed in Table \ref{tab:ex2_parameters}. 

Two DR requests are considered for each day, with random start and duration, so that $\underline{t}_j$ and $\overline{t}_j$ can differ each day ($\underline{t}_2>\overline{t}_1$). 
The lower and upper energy bounds for DR requests are fixed, and they equal to $\underline{E}_1^{DR}=\underline{E}_2^{DR}=-10000~\textnormal{kWh}$, $\overline{E}_1^{DR}=10000~\textnormal{kWh}$ and $\overline{E}_2^{DR}=50000~\textnormal{kWh}$, while both $\overline{\gamma}_1$ and $\overline{\gamma}_2$ are fixed to 3000\euro. For a representative day (day 3), the total load of the community and the total generation by producers not equipped with a BESS are reported in Fig. \ref{fig:ex2_load}, while the net energy injected into the grid is depicted in Fig. \ref{fig:ex2_NET}. Table \ref{tab:ex2_Results} summarizes the profit of each unit both when working individually and when joining the community, for three days of simulation. Finally, the overall extra profit and DR reward of all entities are reported in Table \ref{tab:ex2_amounts}.

{\small 
\begin{table*}[t!]
\renewcommand{\arraystretch}{0.8} % Adjusts the height of all rows
    \centering
    % Table for Entities 1 to 10 (Transposed)
    \begin{tabular}{|c|c|c|c|c|c|c|c|c|c|c|}
    \hline
    Entity & 1 & 2 & 3 & 4 & 5 & 6 & 7 & 8 & 9 & 10 \\
    \hline
    $\overline{P}_u$ & 548 & 412 & 652 & 592 & 364 & 320 & 468 & 588 & 540 & 556 \\
    \hline
    $\overline{S}_u$ & 603 & 277 & 744 & 377 & 291 & 169 & 344 & 711 & 298 & 717 \\
    \hline 
    \hline
    Entity & 11 & 12 & 13 & 14 & 15 & 16 & 17 & 18 & 19 & 20 \\
    \hline
    $\overline{P}_u$ & 516 & 600 & 412 & 400 & 356 & 360 & 472 & 416& 476 & 464 \\
    \hline
    $\overline{S}_u$ & 403 & 563 & 334 & 342 & 370 & 377 & 324 & 245 & 270 & 253 \\
    \hline
    \hline
    Entity & 21 & 22 & 23 & 24 & 25 & 26 & 27 & 28 & 29 & 30 \\
    \hline
    $\overline{P}_u$ & 524 & 696 & 584 & 440 & 328 & 528& 320 & 364 & 648 & 388 \\
    \hline
    $\overline{S}_u$ & 474 & 846 & 624 & 344 & 213 & 527 & 181 & 435 & 499 & 402 \\
    \hline
    \end{tabular}
    \caption{Example 2. Values of entity parameters.}
    \label{tab:ex2_parameters}
\end{table*}}

}
\end{example}

% \begin{table}[t!]
%     \centering
%     \begin{tabular}{|c|c|c|c|c|c|c|c|c|}
%         \hline
%      {Date} & $\underline t_1$& $\overline t_1$ & $\underline t_2$ &  $\overline t_2$&$\underline{E}_1^{DR}$&$\underline{E}_2^{DR}$&$\overline{E}_1^{DR}$&$\overline{E}_2^{DR}$ \\
%         \hline
%        04-01 & 09:45 & 10:30 & 17:45 & 19:00&-10000&-10000&10000& 50000\\
%          04-02 & 09:15 & 10:45 & 16:15 & 17:30&-10000&-10000&10000& 50000 \\
%          04-03 & 08:45 & 09:30 & 17:45 & 18:30&-10000&-10000&10000& 50000 \\
%          04-04 & 08:30 & 09:15 & 16:15 & 16:45 &-10000&-10000&10000& 50000\\
%          04-05 & 10:30 & 11:30 & 16:00 & 17:00&-10000&-10000&10000& 50000 \\
%         \hline
%     \end{tabular}
%    \caption{Example 2. Set of DR requests in 5 simulation days.}
%     \label{tab:ex2_DR}
% \end{table}
\begin{table*}[t!]
    \centering
    \small % Reduce the font size of the table
    \renewcommand{\arraystretch}{1.2} % Adjusts the height of all rows
    \setlength{\tabcolsep}{2pt} % Adjust column separation
    \begin{tabular}{|c|c|c|c|c|c|c||c|c|c|c|c|c|c|}
    \cline{1-13} 
       {Date} & {$J^*_{1,0}$} & {$J^*_{2,0}$} & {$J^*_{3,0}$} & {$J^*_{4,0}$} & {$J^*_{5,0}$} & {$J^*_{6,0}$} & $J^*_{1}$ & $J^*_{2}$ & $J^*_{3}$ & $J^*_{4}$ & $J^*_{5}$ & $J^*_{6}$ \\
    \hline
    04-01 & 242.18 & 180.15 & 288.2 & 258.55 & 159.79 & 139.27 & 351.66 & 261.59 & 418.48 & 375.43 & 232.02 & 202.24 \\
    \hline
    04-02 & 194.64 & 146.24 & 231.58 & 210.07 & 129.28 & 113.1 & 281.56 & 211.55 & 334.99 & 303.88 & 187.02 & 163.62 \\
    \hline
    04-03 & 222.17 & 165.89 & 264.33 & 238.06 & 147.19 & 128.2 & 298.25 & 222.69 & 354.85 & 319.57 & 197.59 & 172.1 \\
    \hline
    \hline
    \cline{1-13} 
       {Date} & {$J^*_{7,0}$} & {$J^*_{8,0}$} & {$J^*_{9,0}$} & {$J^*_{10,0}$} & {$J^*_{11,0}$} & {$J^*_{12,0}$} & $J^*_{7}$ & $J^*_{8}$ & $J^*_{9}$ & $J^*_{10}$ & $J^*_{11}$ & $J^*_{12}$ \\
    \hline
    04-01 & 205.04 & 259.91 & 235.2 & 245.76 & 226.41 & 264.58 & 297.73 & 377.40 & 341.52 & 356.86 & 328.76 & 384.19 \\
    \hline
    04-02 & 166.2 & 208.84 & 191.04 & 197.48 & 183.27 & 213.10 & 240.43 & 302.11 & 276.35 & 285.67 & 265.11 & 308.27 \\
    \hline
    04-03 & 188.84 & 238.39 & 216.5 & 225.41 & 208.55 & 243.2 & 253.50 & 320.02 & 290.64 & 302.60 & 279.96 & 326.48 \\
    \hline
    \hline
    \cline{1-13} 
       {Date} & {$J^*_{13,0}$} & {$J^*_{14,0}$} & {$J^*_{15,0}$} & {$J^*_{16,0}$} & {$J^*_{17,0}$} & {$J^*_{18,0}$} & $J^*_{13}$ & $J^*_{14}$ & $J^*_{15}$ & $J^*_{16}$ & $J^*_{17}$ & $J^*_{18}$ \\
    \hline
    04-01 & 180.94 & 175.92 & 157.26 & 159.04 & 206.48 & 181.41 & 262.74 & 255.44 & 228.36 & 230.94 & 299.83 & 263.42 \\
    \hline
    04-02 & 146.33 & 142.07 & 126.44 & 127.86 & 167.56 & 147.39 & 211.68 & 205.51 & 182.91 & 184.96 & 242.39 & 213.21 \\
    \hline
    04-03 & 166.68 & 162.04 & 144.33 & 145.95 & 190.14 & 167.01 & 223.75 & 217.52 & 193.75 & 195.93 & 255.25 & 224.20 \\
    \hline
    \hline
    \cline{1-13} 
       {Date} & {$J^*_{19,0}$} & {$J^*_{20,0}$} & {$J^*_{21,0}$} & {$J^*_{22,0}$} & {$J^*_{23,0}$} & {$J^*_{24,0}$} & $J^*_{19}$ & $J^*_{20}$ & $J^*_{21}$ & $J^*_{22}$ & $J^*_{23}$ & $J^*_{24}$ \\
    \hline       
    04-01 & 207.44 & 202.06 & 230.81 & 307.65 & 258.03 & 193.06 & 301.21 & 293.40 & 335.16 & 446.72 & 374.68 & 280.35 \\
    \hline
    04-02 & 168.51 & 164.12 & 186.11 & 247.2 & 207.42 & 156.27 & 243.76 & 237.40 & 269.22 & 357.60 & 300.05 & 226.07 \\
    \hline
    04-03 & 190.96 & 186.0 & 212.34 & 282.17 & 236.76 & 177.84 & 256.35 & 249.69 & 285.06 & 378.79 & 317.84 & 238.73 \\
    \hline
    \hline
    \cline{1-13} 
       {Date} & {$J^*_{25,0}$} & {$J^*_{26,0}$} & {$J^*_{27,0}$} & {$J^*_{28,0}$} & {$J^*_{29,0}$} & {$J^*_{30,0}$} & $J^*_{25}$ & $J^*_{26}$ & $J^*_{27}$ & $J^*_{28}$ & $J^*_{29}$ & $J^*_{30}$ \\
    \hline
    04-01 & 143.31 & 233.18 & 139.44 & 160.9 & 284.22 & 171.39 & 208.10 & 338.60 & 202.48 & 233.63 & 412.71 & 248.87 \\
    \hline
    04-02 & 116.4 & 187.53 & 113.28 & 129.28 & 230.16 & 137.80 & 168.39 & 271.28 & 163.86 & 187.02 & 332.93 & 199.35 \\
    \hline
    04-03 & 131.96 & 214.06 & 128.37 & 147.57 & 261.79 & 157.30 & 177.15 & 287.36 & 172.32 & 198.10 & 351.43 & 211.17 \\
    \hline
    \end{tabular}
    \caption{Example 2. Profit of all entities acting individually ($J^*_{i,0}$) and joining the community ($J^*_i$), $i=1,\ldots,30$ under $H^E$, in 3 simulation days.}
    \label{tab:ex2_Results}
\end{table*}

\begin{table*}[ht!]
    \centering
    \centering
     \small % Reduce the font size of the table
    \renewcommand{\arraystretch}{1.2} % Adjusts the height of all rows
    \setlength{\tabcolsep}{2pt} % Adjust column separation
    \begin{tabular}{||c||c||c|c||c|c||c|c|c||}
        \cline{3-6}  
        \multicolumn{2}{c|}{}&
    \multicolumn{2}{c||}{$ H = H^{E}$} & \multicolumn{2}{c||}{ $H = H^{M} $} \\
        \cline{1-6}
           {Date} & $J^*_{0}$ &   $\displaystyle \sum_{u \in \mathcal{U}}\delta_u$&  $\ds\sum_{j:\mathcal{R}_j\in\cR}\gamma^*_j$ &$\displaystyle\sum_{u \in \mathcal{U}}\delta_u$&  $\ds\sum_{j:\mathcal{R}_j\in\cR}\gamma^*_j$  \\
        \hline
        04-01 &   
6297.61& 2846.95 & 3996.31& 2550.03& 3996.31\\
        \hline
        04-02 & 5086.64 & 2271.51 & 2469.96 & 1937.66&2853.52\\
        \hline
        04-03 & 5790.00 & 1982.63 & 2573.94 &1767.63 & 2601.74\\
        \hline
        04-04 &  5829.33 & 1884.17 & 2415.16 & 1559.73& 2525.93\\
        \hline
        04-05 &  5378.23 & 3006.71 & 3118.36& 2667.43& 3488.99\\
        \hline
    \end{tabular}
    \caption{Example 2. Total daily amounts related to two different objective functions $H$, in 5 simulated days.}
    \label{tab:ex2_amounts}
\end{table*}


\section{Discussion} \label{sec:discussion}

In this section, the results of the simulations provided in the above examples are commented in detail. First, let us focus on Example 1.
From Table \ref{tab:table2}, one can notice that the total extra profit of entities ($\delta_1+\delta_2$) is consistently higher under the objective function $H^E$ compared to $H^M$, across all days. On the contrary, the community daily DR reward ($\gamma_1^*+\gamma_2^*$), and hence the revenue of the REC manager, is higher when the objective function $H^M$ is employed. Such results are consistent with the goals of the two objective functions, as reported in Section \ref{sec:Performanceindices}.
% From Table \ref{tab:table2}, it is apparent that the total extra profit of entities is consistently higher under the objective function $H^E$ compared to $H^M$, across all days. 
% Such behavior is reversed for the total daily DR reward at REC level, where $H^M$ provides higher values than $H^E$. This outcome is reasonable because the goal of the objective function $H^M$ is to maximize the revenue of the REC manager which corresponds to maximizing the total DR reward.
As shown in Table \ref{tab:table3}, the proposed DR reward redistribution approach ensures that all entities receive a reward proportional to the individual profits obtained when operating independently. However, the factor $\rho$ is considerably smaller when considering the objective function $H^M$, which implies a reduction of the extra profit for all entities.

Fig. \ref{fig:ex1_alone}, which refers to standalone operation, reveals that during the early hours of the day, when the selling price is low, the entity prioritizes storing its energy production, while discharging starts occurring around 10:00 to take advantage higher selling prices. By 14:00, the BESS is fully discharged to avoid selling the stored energy at lower prices. On the other hand, when joining the REC, both entities stop charging and start exporting their energy production once the first DR request horizon begins (see Fig. \ref{fig:ex1_REC}). During this period, Entity 1 exports 294 kWh, while Entity 2 exports 126 kWh. After the first DR request ends, entities resume charging due to low energy price and later export to the grid during the period with highest energy price. Between 15:15 and 17:00 both entities start charging again, so that their storage systems can be discharged during the second DR request, as expected. During the second DR request, this strategy enables Entities 1 and 2 to provide 729 kWh and 318 kWh, respectively.

From Fig. \ref{fig:ex1_delta} it is apparent that both ${\delta}_1$ and ${\delta}_2$ are higher when the objective function is set to $H^E$ compared to $H^M$. This result is expected, as $H^E$ is designed to maximize the total revenue for the entities, whereas $H^M$ focuses on maximizing the revenue of the REC manager, which naturally leads to lower extra profits for the entities. Fig. \ref{fig:ex1_delta} also demonstrates a fair redistribution of rewards among the entities, such that the additional profit be proportional to their respective individual revenuse when operating independently, i.e., outside the REC. This ensures that the additional profits $\delta_1$ and $\delta_2$ are allocated equitably, reflecting the contribution provided by each entity to the REC. Notice that, when the objective function $H^M$ is used, there are some days such that $\delta_1=\delta_2=0$, that is, the profit of the entities joining the REC equals the profit they would gain if acting individually outside the REC.

Regarding the DR rewards assigned to the entities, Fig. \ref{fig:ex1_DR} shows that when optimizing $H^M$ the total community reward is no less than that obtained by using $H^E$, as expected. However, even if optimizing $H^M$, the REC is in general unable to provide the upper DR request energy bounds (i.e., $\gamma_1^*+\gamma_2^* < \overline\gamma_1+\overline\gamma_2$). Clearly, the sum of the rewards assigned to the two entities is less than the total reward received by the community, since a fraction of it is retained by the REC manager according to \eqref{eq:DR_alpha}.

The sensitivity of the proposed method with respect to the DR reward bound $\overline{\gamma}_j$ is explored in Table \ref{tab:4}. When the DR reward bound $\overline{\gamma}_j$ is increased from 40 to 100\euro, the achieved DR rewards $\gamma_1^*+\gamma_2^*$ under both objectives $H^{E}$ and $H^{M}$ are the same. This basically means that the DR reward is large enough to make fulfillment of DR requests always advantageous regardless of all costs and energy losses arising when operating storage. So, any further increase in $\overline{\gamma}_j$ will provide the same BESS control commands and consequently the same amount of energy injected into the grid.

Concerning Example 2, in Fig. \ref{fig:ex2_NET} one can observe that during the DR periods entities discharge their storage systems to increase the injected energy into the grid. Thanks to this operation, the REC receives a monetary reward which can be shared among entities, allowing them to substantially increase their profit compared to their baseline, see Tables \ref{tab:ex2_Results} and \ref{tab:ex2_amounts}. Regarding the two considered objective functions $H^E$ and $H^M$, they yield behaviors similar to those in Example 1, favouring the entity profit and the total DR reward, respectively.

Example 2 allows one to evaluate the computational burden of the proposed procedure, showing that the algorithm is computationally tractable even in presence of a large number of entities. In fact, the average time needed to solve 30 instances of Problem \ref{pb:alone} and one instance of Problem \ref{pb:manager_new} for a day is about 0.15 seconds on average, allowing this technique to be practically adopted in real-world scenarios. The low computational effort of the whole procedure is due to two reasons: first, Problem \ref{pb:alone} which needs to be solved for each entity is a linear program, and so it can be efficiently solved by standard tools; second, although Problem \ref{pb:manager_new} is a MILP, the low number of involved integer variables allows it to be quickly solved at the optimum. 

\section{Conclusion and future research}\label{sec:conclusions}
In this paper, the potential of coordinating storage operations inside a REC in the presence of DR programs has been investigated. Under the assumption that the REC is involved in price-volume DR programs, a novel 3-step procedure has been proposed to optimize individual storage operation with respect to objective functions that represent overall community benefit. The proposed approach guarantees both an increased profit for REC producers compared to optimally acting outside the community, and the redistribution of DR rewards among participants according to a fairness principle.
The optimization procedure involves the solution of an LP for each schedulable entity and of one MILP with only few integer variables irrespective of the community size, thus making the approach viable for large-scale problem instances. Extensive numerical simulations have been provided showing the effectiveness of the proposed method and comparing the results under two different objective functions.

Future work will focus on the extension of the proposed procedure to different performance indices and reward redistribution policies. Other developments may address the presence of REC prosumers and stochastic behavior of participating entities, as well as uncertainty affecting load and generation profiles.

\section*{CRediT authorship contribution statement}
\textbf{Gianni Bianchini:} Conceptualization, Investigation, Methodology, Writing – original draft, Writing – review \& editing.
\textbf{Marco Casini:} Conceptualization, Funding acquisition, Investigation, Methodology, Writing – original draft, Writing – review \& editing.
\textbf{Milad Gholami:} Data curation, Investigation, Software, Validation, Writing – original draft.

\section*{Declaration of Competing Interest}
The authors declare that they have no known competing financial interests or personal relationships that could have appeared to influence the work reported in this paper.

\section*{Data availability}
Data will be made available on request.

\section*{Acknowledgements}
The project has been funded by the EUROPEAN UNION - NEXT GENERATION EU, MISSION 4, COMPONENT 1, CUP J53D23000760006 (Title project: SHESS4REC - Smart Hybrid Energy Storage Systems for Renewable Energy Communities).

\end{document}